\documentclass[aps,preprint,a4paper,floatfix]{revtex4-1}

\usepackage{graphicx,color}
\usepackage{amsmath,amssymb,amsfonts,amsthm,amscd,bm}
\usepackage{booktabs}
\usepackage{cancel}

\usepackage{amsmath}
\usepackage{amssymb}
\usepackage{amsfonts}
\usepackage{amsthm}
\usepackage{amscd}
\usepackage{bm}
\usepackage{graphicx}
\usepackage{booktabs}
\usepackage{listings}
\usepackage{xcolor,colortbl}
\usepackage{float}

\def\*{\star}
\def\[{\left[}
\def\]{\right]}
\def\({\left(}      

\def\){\right)}

\def\frac#1#2{\dfrac{#1}{#2}}
\def\inv#1{\dfrac{1}{#1}}

\def\2pi{\hbox{$2\pi i$}}

\def\dsl{\raise.15ex\hbox{/}\kern-.57em\partial}
\def\Dsl{\,\raise.15ex\hbox{/}\mkern-.13.5mu D}
\def\2pi{\hbox{$2\pi i$}}

\def\dsl{\raise.15ex\hbox{/}\kern-.57em\partial}
\def\Dsl{\,\raise.15ex\hbox{/}\mkern-.13.5mu D}
\font\numbers=cmss12
\font\upright=cmu10 scaled\magstep1
\def\stroke{\vrule height8pt width0.4pt depth-0.1pt}
\def\topfleck{\vrule height8pt width0.5pt depth-5.9pt}
\def\botfleck{\vrule height2pt width0.5pt depth0.1pt}
\def\Zmath{\vcenter{\hbox{\numbers\rlap{\rlap{Z}\kern
    0.8pt\topfleck}\kern 2.2pt
    \rlap Z\kern 6pt\botfleck\kern 1pt}}}
\def\Qmath{
    \vcenter{\hbox{\upright\rlap{\rlap{Q}\kern3.8pt\stroke}\phantom{Q}}}}
\def\Nmath{\vcenter{\hbox{\upright\rlap{I}\kern 1.7pt N}}}
\def\Cmath{\vcenter{\hbox{\upright\rlap{\rlap{C}\kern
                   3.8pt\stroke}\phantom{C}}}}
\def\Rmath{\vcenter{\hbox{\upright\rlap{I}\kern 1.7pt R}}}
\def\Z{\ifmmode\Zmath\else$\Zmath$\fi}
\def\Q{\ifmmode\Qmath\else$\Qmath$\fi}
\def\N{\ifmmode\Nmath\else$\Nmath$\fi}
\def\C{\ifmmode\Cmath\else$\Cmath$\fi}
\def\R{\ifmmode\Rmath\else$\Rmath$\fi}
\def\barray{\begin{eqnarray}}
\def\earray{\end{eqnarray}}
\def\beq{\begin{equation}}
\def\eeq{\end{equation}}

\def\AA{\leavevmode\setbox0=\hbox{h}
\dimen0=\ht0 \advance\dimen0 by-1ex\rlap{\raise.67\dimen0\hbox{\char'27}}A}

\def\Arg{{\rm Arg}\,}

\def\kappa{c}

\makeatletter
\g@addto@macro\bfseries{\boldmath}
\makeatother

\makeatletter
\renewcommand*\env@matrix[1][\arraystretch]{%
  \edef\arraystretch{#1}%
  \hskip -\arraycolsep
  \let\@ifnextchar\new@ifnextchar
  \array{*\c@MaxMatrixCols c}}
\makeatother

\theoremstyle{plain}

\theoremstyle{remark}
\newtheorem{remark}{Remark}

\newtheorem{hypothesis}{Hypothesis}

\def\Q{\mathbb{Q}}
\def\Z{\mathbb{Z}}
\def\N{\mathbb{N}}

\def\limsup{{\underset{n\to \infty} {\rm lim ~ sup}}}
\def\limsupno{{\rm lim~sup}}

\begin{document}

\title{On large gaps between zeros of $L$-functions from  branches}

\author{
 Andr\'e  LeClair\footnote{andre.leclair@gmail.com}
}
\affiliation{Cornell University, Physics Department, Ithaca, NY 14850} 

\bigskip

\begin{abstract}

It is commonly believed that the normalized gaps between consecutive ordinates $t_n$ of the zeros of the Riemann zeta function
on the critical line can be arbitrarily large. In particular, drawing on analogies with random matrix theory,  it has been conjectured that 
$$\lambda' =\limsup  ~\( t_{n+1} - t_n \) \frac{ \log( t_n /2 \pi e)}{2\pi}$$ 
equals $\infty$.   
In this article we provide  arguments,  although not a rigorous proof,
that $\lambda'$ is finite.    Conditional on  the Riemann Hypothesis,  we show that if there are no  changes of branch between consecutive zeros
then $\lambda' \leq 3$,  otherwise $\lambda'$ is allowed to be greater than $3$.   
 Additional arguments lead us to propose  $\lambda'\leq 5$.
    This proposal  is consistent with numerous
calculations  that place  lower   bounds on $\lambda'$.    
We present the generalization of this result to all Dirichlet $L$-functions and those based on cusp forms.   

\end{abstract}

\maketitle

Let $\rho = \tfrac{1}{2} + i t_n$ denote the $n$-th zero of the Riemann zeta function
on the upper critical line,
with $n=1,2,\ldots$.     An interesting question which has received a great deal of attention
concerns just how large the gaps between consecutive $t_n$  can be.    Since 
on average  the gaps are known to be $2 \pi/\log t_n$,   one is led to study the appropriately normalized gaps  
\beq
\label{gn}
g_n =  \( t_{n+1} - t_n \) \frac{\log t_n}{2 \pi} 
\eeq
and 
\beq
\label{lambda}
\lambda = \limsup ~ g_n
\eeq
Montgomery conjectured that $\lambda = \infty$ based on analogies  with random matrix
theory \cite{Montgomery,OdlyzkoPair}.    

Obtaining {\it lower} bounds on $\lambda$ is a very difficult problem.  
Extensive analysis by many authors indicates that if $\lambda$ is indeed $\infty$,  then it must
approach it very slowly,  if at all.     
The first unconditional result is due to Hall \cite{Hall}:  $\lambda>2.34$.   Assuming the
Riemann Hypothesis (RH),  
Montgomery and Odlyzko  \cite{MO} obtained $\lambda > 1.9799$,  and Conrey,  Ghosh,
and Gonek  \cite{CGG1} improved the result to $\lambda > 2.337$.   If one assumes 
the generalized Riemann Hypothesis (GRH),  one can do slightly better.   The previously mentioned authors obtained 
$\lambda > 2.68$ \cite{CGG2}.    
Ng \cite{Ng} obtained $\lambda > 3$ and Bui \cite{Bui} showed $\lambda > 3.033$.   
The current best unconditional result is due to Bui and Milinovich, $\lambda > 3.18$,  based on the method of Hall \cite{Milinovich}.  
All of these results were obtained by similar methods  involving  studying higher moments of the 
zeta function.     Despite these difficult  analyses,  the results have improved only incrementally,  
and are still very far from the expected $\lambda = \infty$.
   It should be mentioned that Hughes \cite{Hughes}  has proposed some ideas
on how to perhaps get to $\lambda = \infty$.  

 The above  results indicate that either new methods are needed,
or the conjecture $\lambda = \infty$ is false.    In this short note,  we will argue the latter
from a simple argument based on the branches of the zeta function.    Thus the main goal of this article is to
propose an {\it upper}  bound on $\lambda$.  

Let us modify slightly the definition of $g_n$ and $\lambda$:
\beq
\label{lambdap}
g'_n = \( t_{n+1} - t_n \) \frac{ \log( t_n /2 \pi e)}{2\pi}, ~~~~~ \lambda' = \limsup ~ g'_n
\eeq
The extra $2\pi e$  is inconsequential in the limit $n \to \infty$ where $\lambda'=\lambda$,  but as we will explain,  
it is  more instructive  for performing numerical checks of our ideas  at large but finite $n$.   Also,  integer values of $\lambda'$ will  have
a special significance.    
 As we will show, assuming the Riemann Hypothesis,    a simple argument leads to $\lambda' \leq 3$ if there are no changes of branch
 between  consecutive zeros.   The value $3$  is already very close  to the 
 best lower bound $\lambda'> 3.18$ \cite{Milinovich}.      
Additional  arguments lead to the  much more specific  proposal 
\beq
 \lambda' \leq  5
\eeq
To our knowledge, this is  the first proposal for an {\it upper} bound  on  the normalized gaps. 
Although we are unable to  provide a  rigorous proof,    it is worthwhile to elucidate   the arguments leading to this proposal,   
since they are new and lead to a definite but  unexpected result which is in the opposite direction of previous results.  
  Furthermore,  it appears  to be  closer to the 
reality suggested by numerical results.   We have checked numerically  that $g'_n < 3$   for all $n<10^9$. 
For clarity,  we itemize the hypotheses behind our proposal which we could not rigorously prove.   

Below we will present the generalization  of the  above result to all $L$-functions based on Dirichlet characters and  
cusp  (modular) forms.    In order to avoid introducing too many definitions that are well-known for these cases,  we will present the main arguments for 
 the Riemann zeta function itself, 
and subsequently  simply state  its  generalization to these other $L$-functions.    As will be clear,  our arguments assume the RH appropriate to the $L$-function in question is true,
although indirectly.   For instance,  we do not need to assume the GRH  to study $\lambda'$ for Riemann zeta itself,  but only the original RH.

Let $\vartheta (t)$ denote the Riemann-Siegel $\vartheta$ function:
\beq
\label{RStheta}
\vartheta (t) =  \Im \log \Gamma(\tfrac{1}{4} + \tfrac{it}{2} ) - t \log \sqrt{\pi} 
\eeq
and $a(t)$ the argument of the zeta function defined in the following specific manner:
 \beq
 \label{aoft}
 a(t) = \lim_{\delta \to 0^+} \arg \zeta (\tfrac{1}{2} + \delta + it ). 
 \eeq
In \cite{electrostatic,Trans} it was proposed that the  ordinate $t_n$ of the $n$-th zeta zero on the upper critical line,
with $n=1,2,3\ldots$, satisfies the exact  transcendental equation 
\beq
\label{trans}
\vartheta (t_n) + a(t_n) = (n - \tfrac{3}{2})\pi
\eeq
Several remarks regarding this equation are in order.   The equation was obtained by putting the zeros in one to one correspondence with the
zeros of the  cosine function.     In obtaining this formula,  we did not assume  Backlund's counting formula for $N(T)$,  i.e. the number of zeros up to height $T$ in the entire critical strip,
$N(T) = \vartheta (T)/\pi + 1 + S(T)$,  where $S(T) = \inv{\pi} \arg \zeta ( \tfrac{1}{2} + iT )$.    The equation \eqref{trans} contains
more information that $N(T)$.   It is important to note that $S(T)$ and $a(t)$ are not  equivalent.   For instance
$S(T)$ is not defined at  the ordinate of a zero,  unlike $a(t)$. See \cite{SofTcite} for a review of known properties of $S(T)$.   It is also important that $a(t)$ is defined precisely as in \eqref{aoft},  i.e. 
as a limit from the right of the critical line with $\delta$ not allowed to be strictly zero.    The behavior of $a(t)$ would be completely different
if it were defined as a limit from the left of the line or along the line.  (See Remark \ref{remark1}.)        
 If there is a unique solution to \eqref{trans} for every positive integer $n$,   then the RH is true since  this implies that the number of zeros on the critical line
 saturate the formula $N(T)$ \cite{electrostatic,Trans}.       If one ignores the $a(t_n)$ term there is in fact a unique solution for every $n$ that 
 can be expressed in terms of the Lambert $W$ function to a very good approximation.    
However we  were unable to rigorously prove there exists a unique solution including the $a(t_n)$ term;  thus  we take it as our first hypothesis:
 
 \bigskip
 \begin{hypothesis} 
 \label{hyp1}
 There is a unique solution to the equation \eqref{trans} for every positive integer $n$.     This effectively assumes the RH and that the zeros are simple\cite{Trans}.   
 \end{hypothesis}
 \bigskip

 Using the Stirling formula, 
 \beq
\label{RSasym}
\vartheta (t)  = \frac{t}{2} \log \( \frac{t}{2 \pi e} \) -\frac{\pi}{8} + O(1/t).
\eeq
Henceforth it is implicit that we are considering the limit of large $n$ since we use the above
asymptotic behavior,  i.e. we will not always display the $\limsupno$.     Assuming Hypothesis \ref{hyp1},
we can consider the difference of the equation for two consecutive $n$'s.  
Clearly one has $(t_{n+1} - t_n )\log t_n < t_{n+1} \log t_{n+1} - t_n \log t_n $,  thus
\beq
\label{gpb1} 
 g'_n  < 1 - \inv{\pi} \( a(t_{n+1}) - a (t_n) \) 
 \eeq

\def\dbmax{\Delta b_{\rm max}}

First of all,  the equation \eqref{gpb1} shows that 
if $a(t)$ is finite,  then so is $\lambda'$.   See Remark \ref{remark1} for comments 
concerning this finiteness.  
 In order to  constrain $\lambda'$,  we need one additional hypothesis.  
 Let us write 
\beq
\label{branch}
a(t_n) = A_n+ 2\pi  b_n, ~~~~~b_n \in \Zmath
\eeq
where $A_n$ denotes  the principal branch:
$A_n = \lim_{\delta \to 0^+} \Arg \zeta (\tfrac{1}{2} + \delta + i t_n)$
with $-\pi < A_n < \pi$.   We will refer to $b_n$ as the branch number,  where
$b_n =0$ corresponds to the principal branch.     
Since $|A_{n+1} - A_n| \leq 2 \pi$, equation \eqref{gpb1} implies 
\beq
\label{gpb2} 
g'_n \leq 3 + \dbmax
\eeq
where 
\beq
\label{bmax} 
\dbmax = {\rm max}  ~  (b_n - b_{n+1} )
\eeq
From this equation we deduce that if there are no changes of branch between consecutive zeros,  then
$g'_n \leq 3$ and $\lambda' \leq 3$.    Changes of branch allow $\lambda' \geq 3$,  but do not necessarily imply it.    Based on the result 
$\lambda > 3.18$  in \cite{Milinovich}  we should conclude that branch changes do occur;  however all the way up 
to $n=10^9$ we found $g'_n < 3$ (see below).     The main shortcoming of our approach is that is that it leads to  bounds on $\lambda'$ 
that are only integers;  it would be impossible to constrain it further to specific values like $\lambda > 3.18$ with this reasoning.  
However we think the approach is useful since it endows values like $\lambda'=3$ some special significance.

 An {\it upper} bound on 
$\lambda'$ would follow from a bound on how much the branch number $b_n$ can decrease between
consecutive zeros.    Let us assume that $a(t)/\pi$ behaves like $S(t)$ as far a branch changes are concerned.   Between zeros,  $a(t)$ always decreases.    Assuming Hypothesis \ref{hyp1},  since the zeros are simple,  the only place where the branch number of
$a(t)/\pi $ can {\it increase}  is at the jumps by $1$  at each zero.    Therefore there is an upper bound to increases of branch number: 
\beq
\label{bup}
b_{n+1} - b_n \leq 1
\eeq
Now since the average of  $g_n'$ is known to be $1$,  the average of $a(t_{n+1}) - a(t_n)$ is zero.    This suggests that  the average of 
$b_{n+1} - b_n$ equals zero.   
A reasonable assumption then is  that  there is no preference between increases verses decreases of branch number.   
This leads to our second hypothesis,  which is likely to be the more difficult one to establish:

\bigskip

\begin{hypothesis}
\label{hyp2} 
There is no dichotomy between increases and decreases of branch number.    Then \eqref{bup} implies  
\beq
|b_{n+1} - b_n | \leq 1
\eeq
\end{hypothesis}

\bigskip

 This hypothesis implies $\dbmax  \leq  1$.  Combined with \eqref{gpb2} this leads to $ \lambda' \leq 5$.   
Note that, in contrast, $\lambda'=\infty$ would require an infinite number of branch changes between two consecutive zeros.   
    
\bigskip

 \begin{remark}
 \label{remark1}
   Based on the validity of the Euler product formula to the right of the critical line (for $L$ functions based on principal characters, the
   product must be truncated in a well-prescribed fashion),  it was 
 argued  in \cite{EPF2,ALzeta} that $a(t)$ is finite,  which is consistent with an upper bound on $\lambda'$:  if $a(t)$ is finite,  then so is $\lambda'$.   
       Let us summarize the arguments that led to this conjecture.    Consider first  $L$-functions
 based on a  non-principal Dirichlet character $\chi$.    We conjectured that $C_N= \sum_{n=1}^N  \chi (p_n)$,  where $p_n$ is the $n$-th prime, 
 behaves like a random walk and is thus $O(\sqrt{N})$ (up to logs).    Further support for this conjecture was obtained using a variant of 
 Cram\'er's random model for the primes \cite{ALCLT}.     Assuming this random walk behavior,  we proved that the Euler product converges 
 to the right of the critical line.   For the present work,   since we are assuming the GRH for non-principal characters,  we can actually use the result 
 $C(x) =\sum_{p\leq x}  \chi (p)  = O(\sqrt{x} \log^2 x)$ (\cite{Davenport} page 125).  The $\log^2 x$ does not spoil the convergence  arguments in \cite{EPF2},  so the GRH implies the validity of the Euler product
 formula to the right of the critical line for non-principal Dirichlet characters.  
 One can therefore  calculate $a(t)$ from the Euler product: 
  \beq
 \label{afinite}
 a(t) = - \lim_{\delta \to 0^+} \Im \sum_p \log \( 1 -  \frac{\chi (p)}{p^{1/2+ \delta + it}} \) 
 \eeq
 If  the product converges,  then $a(t)$ and $\lambda'$  are finite.     
  For principal characters the situation is more subtle and one has to truncate the Euler product
 \cite{Gonek1, Gonek2,ALzeta} at $p=p_{N_c}$ where $N_c \sim t^2$,  however this does not spoil the conclusion that $a(t)$ is finite.   
 In fact,  if one is interested in $\limsupno$,  the truncation is not necessary since $N_c \to \infty$ as $t \to \infty$. 
 Furthermore,  approximating the expression in \eqref{afinite} using the prime number theorem   led us to 
 propose that $a(t)$ is nearly always on the principal branch \cite{ALzeta}  (based on equation (38) in \cite{ALzeta}).   The latter approximation also supports Hypothesis \ref{hyp2},  since there is no  dichotomy between the oscillations above verses below zero.
     Consequently in regions where $a(t_n)$ is on the principal branch, 
 $g'_n < 3$;    we will provide numerical evidence for this below.  
\end{remark}

We can easily extend the above arguments to two additional classes of $L$-functions.    First consider $L$-functions based on 
any primitive Dirichlet character mod $q$.    Denote the  zeros on the critical line as $\rho = \tfrac{1}{2} + i t_n$.    Repeating the above 
arguments using the transcendental equations in \cite{Trans} one finds that the proper normalization depends on $q$:
\beq
\label{lamdaDir}
\lambda' = \limsup 
 \( t_{n+1} - t_n \) \frac{ \log(q \, t_n /2 \pi e)}{2\pi},
\eeq
Next consider $L$ functions based on cusp forms mod $k$.    Here the zeros on the line are $\rho = \tfrac{k}{2} + i t_n$.   
In this case the proper normalization does not depend on $k$ however there is a difference by an overall  factor of $2$:
\beq
\label{cusp}
\lambda' = \limsup 
 \( t_{n+1} - t_n \) \frac{ \log( t_n /2 \pi e)}{\pi},
\eeq
In both cases the previous arguments again imply $ \lambda' \leq 5$.

\begin{figure}[t]
\centering\includegraphics[width=.7\textwidth]{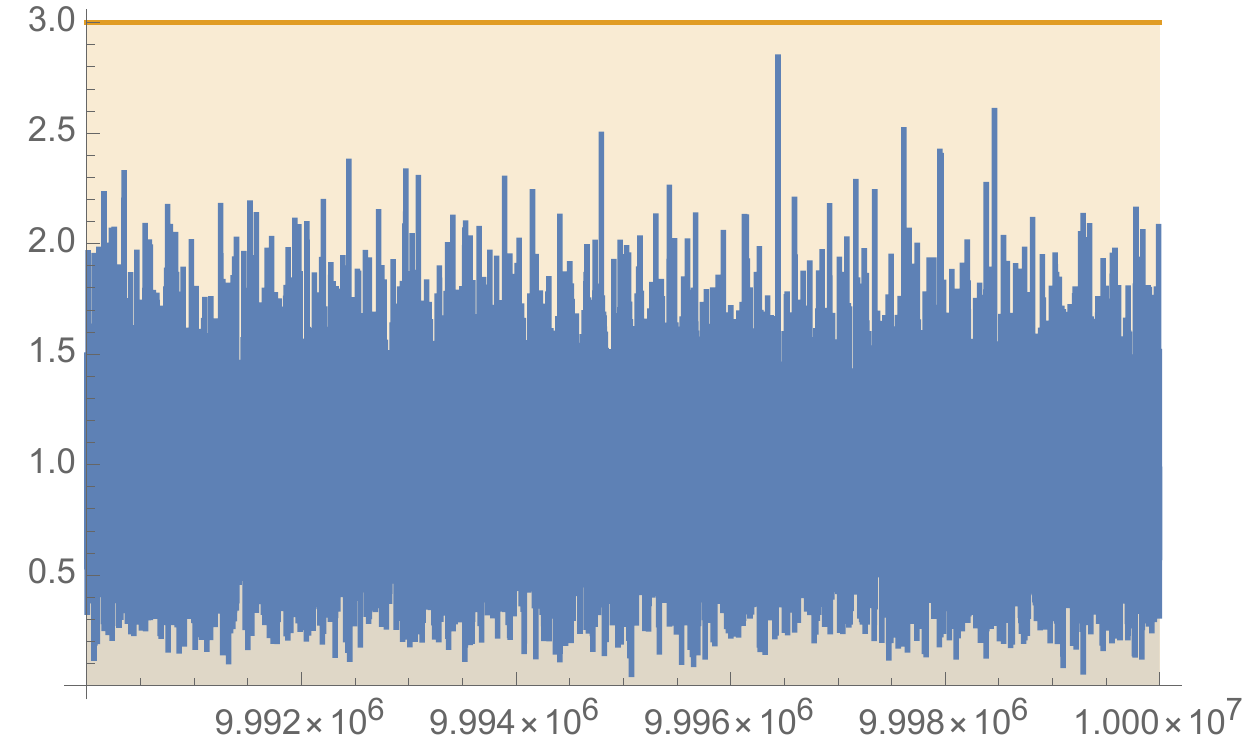}
\caption{The normalized gaps $ g'_n  = \( t_{n+1} - t_n \)  \log( t_n /2 \pi e) /2\pi$  of the zeta function for \\
~~~~~~~~$ 10^7-10,000 < n < 10^7$. }
\label{gaps}
\end{figure}

\begin{figure}[t]
\centering\includegraphics[width=.7\textwidth]{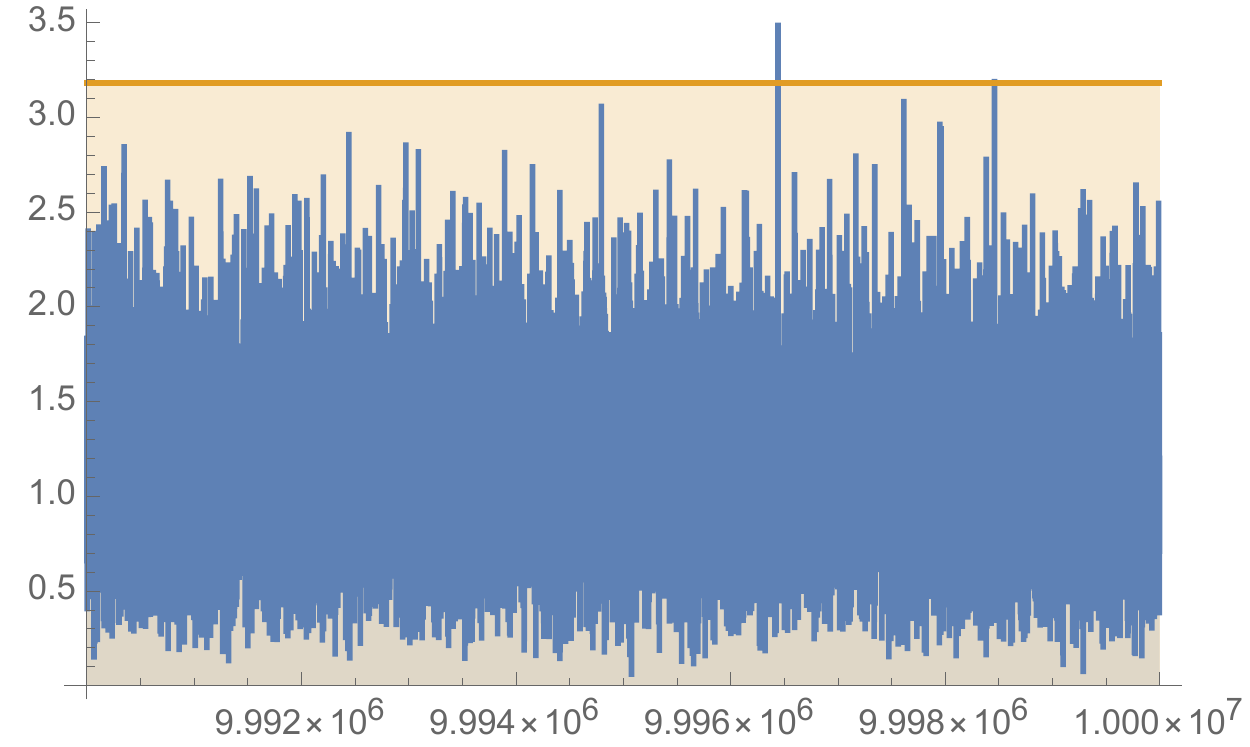}
\caption{The normalized gaps $ g_n  = \( t_{n+1} - t_n \)  \log( t_n ) /2\pi$  of the zeta function for \\
~~~~~~~~$ 10^7-10,000 < n < 10^7$, to be compared with Figure \ref{gaps}.   The horizontal line is  $3.18$,  based on the prediction in \cite{Milinovich}
that $\lambda'>3.18$. }
\label{gapOld}
\end{figure}

Let us now provide numerical support for the above proposals.   We limit ourselves to the zeta function.    
All the way up to $n=10^9$ we found that $g'_n < 3$,  with some cases coming close to this upper bound. 
     In Figure \ref{gaps} we display a region around $n=10^7$.  
The different definitions of $\lambda$ and $\lambda'$ are significant here in that for some of the extreme gaps 
where $g'_n $ is just under $3$,  one can check that  $g_n >3$.    As displayed in  Figure \ref{gapOld},  the rare points where $g_n >3.18$ in this range 
still have $g'_n < 3$.      The explanation of these results is simple based on 
the above ideas:   in this whole range,  $a(t_n)$ is always on the principal branch with $b_n = 0$ for all $n$,  and by \eqref{gpb2}, 
$g'_n < 3$.   (See Remark \ref{remark1}.)   Had we used the definition $g_n$ instead of $g'_n$,  the rare points where $g_n > 3$ would not so easily be explained. 
We also wish to point out that if $a(t_n)$ is never very far from the principal branch,  then $\lambda' \approx 3$,  which suggests that 
Bui and Milinovich's result 
 $\lambda > 3.18$ may actually be close to the  true upper bound.

\section{Acknowledgements} 

We wish to thank Nelson Carella and Micah Milinovich for discussions.


\end{document}